\def\draft{n}
\documentclass{amsart}
\usepackage{fullpage,amssymb,epic,eepic,epsfig,amsbsy,amsmath}

\theoremstyle{plain}

\newtheorem{theorem}{Theorem}
\newtheorem{proposition}{Proposition}[section]
\newtheorem{lemma}[proposition]{Lemma}
\newtheorem{corollary}[proposition]{Corollary}

\newtheorem{conjecture}{Conjecture}

\theoremstyle{definition}
\newtheorem{definition}[proposition]{Definition}

\newtheorem{question}{Question}

\theoremstyle{remark}

\newtheorem{remark}[proposition]{Remark}

\def\printname#1{
	\if\draft y
		\smash{\makebox[0pt]{\hspace{-0.5in}
			\raisebox{8pt}{\tt\tiny #1}}}
	\fi
}

\newcommand{\psdraw}[2]
         {\begin{array}{c} \hspace{-1.3mm}
	\raisebox{-4pt}{\epsfig{figure=draws/#1.eps,width=#2}}
	\hspace{-1.9mm}\end{array}}

\newlength{\standardunitlength}
\catcode`\@=11
\long\def\@makecaption#1#2{%
     \vskip 10pt

\setbox\@tempboxa\hbox{
       \small\sf{\bfcaptionfont #1. }\ignorespaces #2}%
     \ifdim \wd\@tempboxa >\captionwidth {%
         \rightskip=\@captionmargin\leftskip=\@captionmargin
         \unhbox\@tempboxa\par}%
       \else
         \hbox to\hsize{\hfil\box\@tempboxa\hfil}%
     \fi}
\font\bfcaptionfont=cmssbx10 scaled \magstephalf
\newdimen\@captionmargin\@captionmargin=2\parindent
\newdimen\captionwidth\captionwidth=\hsize
\catcode`\@=12

\def\lbl#1{\label{#1}\printname{#1}}

\def\BZ{\mathbb Z}

\def\BQ{\mathbb Q}

\def\BC{\mathbb C}

\def\D{\Delta}

\def\a{\alpha}

\def\s{\sigma}

\def\e{\epsilon}

\def\s{\sigma}

\def\sgn{\operatorname{sgn}}

\def\mat#1#2#3#4{\left(
\begin{matrix}
 #1 & #2  \\
 #3 & #4   
\end{matrix}
\right)}

\def\longto{\longrightarrow}
\def\sig{$\s$ignature}
\def\jump{\mathrm{j}}
\def\jjump{\mathrm{jj}}
\def\mult{\mathrm{mult}}
\def\DivD{\mathrm{Div}_{\D(K)}}
\def\bt{\bar{t}}
\def\pt{\partial}
\def\coeff{\mathrm{coeff}}

\begin{document}


\title[Does the Jones polynomial determine the signature of a knot?]
{Does the Jones polynomial determine the signature of a knot?}

\author{Stavros Garoufalidis}
\address{School of Mathematics \\
         Georgia Institute of Technology \\
         Atlanta, GA 30332-0160, USA \\ 
         {\tt http://www.math.gatech} \newline {\tt .edu/$\sim$stavros } }
\email{stavros@math.gatech.edu}

\thanks{The author was supported in part by National Science Foundation and
Binational Science Foundation. \\
\newline
1991 {\em Mathematics Classification.} Primary 57N10. Secondary 57M25.
\newline
{\em Key words and phrases: Jones polynomial, signature of knots, colored 
Jones function, Alexander polynomial, jump divisor.
}
}

\date{September 28, 2003 \hspace{0.5cm} First edition: September 28, 2003.}


\begin{abstract}
The signature function of a knot is a locally constant integer valued 
function with domain the unit circle. The jumps (i.e., the discontinuities)
of the signature function can occur only
at the roots of the Alexander polynomial on the unit circle. The latter are
important in deforming $U(1)$ representations of knot groups to irreducible
$SU(2)$ representations. Under the assumption that these roots are simple,
we formulate a conjecture that explicitly computes the jumps of the signature 
function in terms of the Jones polynomial of a knot and its parallels.
As evidence, we prove our conjecture for torus knots, and also (using 
computer calculations) for knots with at most $8$ crossings.
We also give a formula for the jump function at simple roots in terms
of relative signs of Alexander polynomials.
\end{abstract}

\maketitle

\tableofcontents


\section{Introduction}
\lbl{sec.intro}

\subsection{The signature function of a knot}
\lbl{sub.history}

A celebrated invariant of a knot $K$ in 3-space is its {\em \sig\ function}
$$
\s(K): S^1 \longto \BZ,
$$ 
defined for complex numbers of absolute value $1$, and taking values in the
set of integers.
The signature function of a knot is a concordance invariant, and plays a 
key role in the study of knots via surgery theory, \cite{L}.

It turns out that the signature function is a locally constant function
away from the (possibly empty) set
$$
\DivD=\{ \rho \in S^1 \, | \D(K)(\rho)=0 \}
$$ 
of roots of the Alexander polynomial on the unit circle. In view of this,
the interesting part of the signature function is its {\em jumping behavior}
on the set $\DivD$.
 
In other words, we may consider the associated {\em jump function} 
$$
\jump(K): \DivD \longto \BZ
$$
defined by $\jump_{\rho_0}(K)=\lim_{\rho \to \rho^+_0}
\s_{\rho}(K)-\lim_{\rho \to \rho^-_0}
\s_{\rho}(K)$.

We may identify the jump function with  a {\em jump divisor}
$\sum_{\rho \in \DivD} \jump_{\rho}(K)[\rho]$ in $S^1$.

Since $1 \not\in \DivD$ and $\s_1(K)=0$, it follows that the jump function
uniquely determines the signature function away from the set $\DivD$.
Since $-1 \not\in \DivD$, it follows in particular that $\jump(K)$ determines
the {\em \sig\ } of the knot $\s_{-1}(K)$.

The signature of a knot may be defined using a Seifert surface of a knot
(see Section \ref{sub.symmetries} below). An intrinsic definition of the
jump function of a knot was given by Milnor \cite{M1,M2}, using the
Blanchfield pairing of the universal abelian cover of a knot. This definition,
among other things, makes evident the role played by the roots of the Alexander
polynomial on the unit circle (as opposed to the rest of
the roots of the Alexander polynomial, which are ignored).

From the point of view of gauge theory and mathematical physics, 
the signature function of a knot may be identified with the
spectral flow of a 1-parameter family of the signature operator, twisted
along abelian (that is, $U(1)$-valued) representations of the knot complement.

The moduli space of $U(1)$ representations of the knot complement is well
understood; it may be identified with the unit circle. On the other hand,
the moduli space of $SU(2)$ representations is less understood, and 
carries 
nontrivial topological information about the knot and its Dehn fillings, as 
was originally discovered by Casson (see \cite{AM}) and also by X-S. Lin;
see \cite{Li}.

One may ask to identify the $U(1)$ representations which {\em deform to} 
irreducible $SU(2)$
representations. Using a linearization argument, Klassen
and Frohman showed that a necessary condition for a $U(1)$ representation
$\rho$ to deform is that $\D(K)(\rho^2)=0$. This brings us to the 
(square of the) set $\DivD$. Conversely, Frohman-Klassen proved sufficiency
provided that the Alexander polynomial has {\em simple roots} on the unit
circle; see \cite{FK}. 
Herald proved sufficiency under the (more relaxed condition that) the jump 
function vanishes nowhere; see \cite{H1,H2}.

It is unknown at present whether sufficiency holds without any further
assumptions.

Let us summarize the two key properties of the jump divisor $\DivD(K)$,
in the spirit of Mazur (see \cite{Ma}):

\begin{itemize}
\item
The jump divisor controls the signature function of a knot.
\item
The jump divisor controls (infinitesimally) deformations of $U(1)$ 
representations of the knot complement to irreducible $SU(2)$ representations.
\end{itemize}

\subsection{The colored Jones function of a knot}
\lbl{sub.coloredJ}

It is a long standing problem to find a formula for the \sig\ function
of a knot in terms of its {\em colored Jones function}. The latter is
a sequence of Jones polynomials associated to a knot. Recall that given
a knot $K$ and a positive integer $n$ (which corresponds to an $n$-dimensional
irreducible representation of $\mathfrak{sl}_2$), one can define a Laurrent
polynomial $J_n(K) \in \BZ[q^{\pm }]$.

In \cite{R2}, Rozansky considered a repackaging of the sequence 
$\{J_n(K)\}$. Namely, he defined a sequence of rational
functions $Q_k(K) \in \BQ(q)$ for $k \geq 0$ with the following properties:
\begin{itemize}
\item
$Q_k(K)=P_k(K)/\D^{2k+1}(K)$ for some polynomials $P_k(K) \in \BZ[q,q^{-1}]$ 
with $P_0(K)=1$ and such that $P_k(K)(q)=P_k(K)(q^{-1})$.
\item
For every $n$ we have:
\begin{equation}
\lbl{eq.euler}
\frac{J_n(K)(q)}{J_n(\text{unknot})(q)}
=\sum_{k=0}^\infty Q_k(q^n) (q-1)^k \in \BQ[[q-1]]
\end{equation}
where $\BQ[[q-1]]$ is the ring of formal power series in $q-1$ with 
rational coefficients 
\end{itemize}

Equation \eqref{eq.euler} is often called the {\em Euler expansion} of the
colored Jones function.
In physical terms, the above expansion is an asymptotic expansion of the 
Chern-Simons
path integral of the knot complement, expanded around a backround
$U(1)$ flat connection.
Thus, philosophically, it should not be a surprise to discover that this
expansion has something to do with the signature of the knot.

For the curious reader, let us point out that Rozansky conjectured such
an expansion for the full Kontsevich integral of a knot, graded by the negative
Euler characteristic of graphs (thus the name, Euler expansion).
This conjecture was proven by Kricker and 
the author; \cite{GK1}. Furthermore, a close relation was discovered between
residues of the rational functions $Q_k$ at roots of unity and the LMO
invariant of cyclic branched coverings of the knot; \cite{GK2}.
In an attempt to understand the Euler expansion, 
a theory of finite type invariants of knots (different from the usual
theory of Vassiliev invariants) was proposed in \cite{GR}. According to that
theory, two knots are $0$-equivalent iff they are $S$-equivalent; \cite{GR}.
Moreover, $Q_k$ is a finite type invariant of type $2k$.

Technically, the Euler expansion of the colored Jones function
is an integrality statement. Namely, it is easy to see that there exist
unique sequence of power series $Q_k(K)(q) \in \BQ[[q-1]]$ for $k \geq 0$
that satisfies Equation \eqref{eq.euler}.
The hard part is to show that these power series are Taylor series expansions
of rational functions with integer coefficients and prescribed denominators.

The statement $P_0(K)=1$ in the leading term
of the Euler expansion is nothing but the Melvin-Morton-Rozansky conjecture,
proven by Bar-Natan and the author in \cite{BG}. Thus, the leading order
term in the Euler expansion is a well-understood topological invariant
of knots. Ever since the Euler expansion was established, it has been
a question to establish a topological understanding of the lower order terms.

\subsection{The conjecture}
\lbl{sub.conjecture}

Consider  $Q(K)(t)=\frac{P(K)(t)}{\D^2(K)(t)} \in \BQ(t)$ where $P(K)=P_1(K)$.
We will think of $Q(K)$ as a function (with singularities) 
defined on the unit circle.

If $\rho=e^{i \theta_0}$ is a root of the Alexander polynomial on $S^1$, we 
may expand $Q(K)(e^{i \theta})$ around $\theta=\theta_0$. The result is a 
power series with lowest term $c_{\rho} (\theta-\theta_0)^{m_{\rho}}$, for some
integer $m_{\rho}$ and some nonzero real number $c_{\rho}$.

\begin{definition}
\lbl{def.jjump}
Let us define the {\em Jones jump function} of a knot $K$
$$
\jjump(K): \DivD \longto \BZ
$$
by 
$$
\jjump_{\rho}(K)= \sgn( c_{\rho}) \, \max\{0, - m_{\rho}\} \, 
\sgn(\mathrm{Im} (\rho))
$$
where $\mathrm{Im}(z) $ is the imaginary part of a complex number $z$ and
$\sgn(x)$ is the {\em sign} of
a real number $x$ is defined by $\sgn(x)=+1,0$ or $-1$ according to 
$x>0, \, x=0$ or $x <0$ respectively.
\end{definition}

\begin{definition}
\lbl{def.simple}
We say that a knot $K$ is {\em simple} if its Alexander polynomial
$\D(K)$ has simple roots on the unit circle.
\end{definition}

\begin{conjecture}
\lbl{conj.1}
If $K$ is simple, then $\jump(K)=\jjump(K)$.
\end{conjecture}

A modest corollary is:

\begin{corollary}
\lbl{cor.1}
If $K$ is simple, Conjecture
\ref{conj.1} implies that the colored Jones function of $K$ determines
the signature $\s_{-1}(K)$.
\end{corollary}

\begin{remark}
\lbl{rem.doesnotmatter}
Notice that $\jump_{\rho}(K)=-\jump_{\bar\rho}(K)$ and 
$\jjump_{\rho}(K)=-\jjump_{\bar\rho}(K)$. Thus, it suffices to check
the conjecture on the upper semicircle.
\end{remark}

\begin{remark}
The conjecture is false if $\D(K)$ has multiple roots (of odd or even
multiplicity). For example,
consider the connected sum $\sharp^n K$ of $n$ right trefoils. Then,
$Q(\sharp^n K)=n Q(K)$ and $\D(\sharp K)=\D(K)^n$.
\end{remark}

We present the following evidence for the conjecture:

\begin{theorem}
\lbl{thm.1}
$\mathrm{(a)}$ Conjecture \ref{conj.1} is true for torus knots, and for
knots with at most 8 crossings. \\
$\mathrm{(b)}$ The Conjecture is compatible with the operations of mirror
image, connected sum (assuming the resulting knot is simple) and $(n,1)$
parallels of knots.
\end{theorem}

En route to establish our results, we give a skein formula that uniquely
characterizes the jump function of simple knots; see Theorem \ref{thm.unique}.

Let us compare Conjecture \ref{conj.1} with existing conjectures about
the structure of the colored Jones function. At the time of the writing,
there are two conjectures that relate the colored Jones function to
hyperbolic geometry. Namely,
\begin{itemize}
\item
The {\em Hyperbolic Volume Conjecture}, 
after Kashaev and J\&J.Murakami,
which states that for a hyperbolic knot $K$,
$$
\lim_{n \to \infty} \frac{ \log| J'_n(K)(e^{2 \pi i/n})|}{n}
=c \, \text{vol}(S^3-K)
$$
where $J'_n(K)=J_n(K)/J_n(\text{unknot})$.
\item
The {\em Characteristic equals deformation variety Conjecture}, due to
the author, which compares the deformation curve of $\mathrm{SL}_2(\BC)$
representations of a knot complement (viewed from the boundary) with a complex
curve which is defined using the recursion relations (with respect to $n$)
of the sequence $\{J_n(K)\}$; see \cite{GL} and \cite{Ga3}.
\end{itemize}

The Hyperbolic Volume Conjecture is an analytic statement, which involves
the existence and identification of a sequence of real numbers.

On the other hand, the Characteristic equals Deformation Variety conjecture
is an algebraic statement, since it is equivalent to the equality of two
polynomials with integer coefficients, one of which is obtained by 
noncommutative elimination, and the other obtained by commutative elimination.

Conjecture \ref{conj.1} appears to be an analytic conjecture, since its basic
ingredients are signs of real numbers. In the field of Quantum Topology, 
analytic conjectures have held the longest.

Let us end the introduction with the following

\begin{question} 
Understand the underlying geometry and perturbative quantum field theory 
behind the Taylor expansion of the $Q$ function (and more generally, Euler
expansion \eqref{eq.euler} of the colored Jones function).
In particular, use the higher order terms $Q_k$ in the expansion 
\eqref{eq.euler} to formulate a conjecture for the jump function of all knots.
\end{question}

\subsection{Acknowledgement}
The author wishes to thank S. Orevkov, L. Rozansky and A. Stoimenov,
and especially J. Levine for helpful conversations.

\section{The signature and the jump function}
\lbl{sec.sj}

\subsection{Symmetries of the jump function}
\lbl{sub.symmetries}

Given a Seifert matrix $V$ of a knot $K$, consider the {\em Hermitian matrix}
$B(t)=(1-t) V + (1-\bar t) V^T$, for $t \in S^1$. The eigenvalues of $B(t)$
are real, and we define $\s_t(K)=\s(B(t))$, where $\s(M)$ denotes the
{\em signature} of a Hermitian matrix $M$. It turns out that $\s(K)$
is independent of the Seifert surface $V$ chosen.
Since $B(t)=(t^{1/2}-t^{-1/2}) A(t)$, where $A(t)=t^{1/2}V -t^{-1/2}V^T$,
and $\det(A(t))=D(K)(t)$ is the {\em symmetrized Alexander polynomial} 
of $K$, it follows that $\s(K)$ is a locally constant function with possible 
jumps along the set $\DivD$. 

The next lemma, which follows from the proof of \cite[Corollary 2]{H1},
summarizes the symmetries of the jump function.

\begin{lemma}
\lbl{lem.1}
If $\rho$ is a root of the Alexander polynomial on $S^1$, then
$|\jump_{\rho}(K)|=2 \, a_{\rho}$, where 
\begin{itemize}
\item
$\a_{\rho}$ is an integer 
\item
$a_{\rho} \leq \mult(\rho,\D(K))$, where $\mult(\rho,\D(K))$ is the
multiplicity of $\rho$ in $\D(K)$, and
\item 
$\a_{\rho} \equiv \mult(\rho,\D(K)) \bmod 2$.
\end{itemize}
Moreover, $\jump_{\rho}(K)=-\jump_{\bar{\rho}}(K)$.
\end{lemma}

In particular, if $K$ is simple, $\jump(K)$ takes values in the set $\{-2,2\}$.
For a precise formula for the jump function in that case, see Theorem 
\ref{thm.jump}.

\subsection{A skein theory for the signature and the jump function}
\lbl{sub.skein}

Let us begin with a useful definition.
A triple of links $(L^+, L^-, L^0)$
is called {\em bordered} if there is an embedded ball $D^3$
in $S^3$ that locally intersects them as in figure \ref{crossing}.

\begin{figure}[htpb]
$$ \psdraw{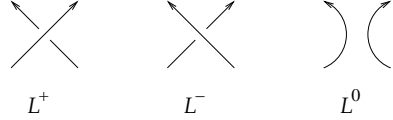}{2in} $$
\caption{A bordered triple of links $(L^+, L^-, L^0)$.}\lbl{crossing}
\end{figure}

If we choose planar projection and a crossing, then a bordered triple 
corresponds to replacing the crossing by a positive, negative or
smoothening. 
Notice that if $L^+$ is a link with $1$ component, 
then $L^-$ and $L^0$ are links with $1$ and $2$ components respectively.

The next lemma computes the change of the signature function with respect
to the change of a crossing, in terms of the sign of the Alexander 
polynomials. 

\begin{lemma}
\lbl{lem.skeins}
If $K$ is a knot, $\rho=e^{i \theta} \in S^1$ such that 
$\D(K^+)(\rho) \D(K^-)(\rho) \neq 0$, then
\begin{equation}
\lbl{eq.skeins}
\s_{\rho}(K^-)-\s_{\rho}(K^+)=
\begin{cases}
2 & \text{if} \quad \D(K^+)(\rho) \D(K^-)(\rho) < 0 \\
0 & \text{if} \quad \D(K^+)(\rho) \D(K^-)(\rho) > 0 .
\end{cases}
\end{equation}
\end{lemma}

\begin{proof}
We can choose Seifert surfaces $(V_+,V_-,V_0)$ for $(K^+,K^-,K^0)$ such
that
$$
V_+= \mat {a}{v_1}{v_2^T}{V_0} \qquad
V_-= \mat {a+1}{v_1}{v_2^T}{V_0}  
$$
where $v_1$ and $v_2$ are some row vectors.
Hermitianizing, we get:
$$
B_+= \mat {2a}{v}{v^T}{B_0} \qquad
B_-= \mat {2a+2-2\cos\theta}{v}{v^T}{B_0}.
$$
Let us call a  triple of Hermitian matrices $(A_+, A_-, A_0)$ 
{\em $\rho$-bordered} if
$$ 
A_+ = \mat {a} {v} {v^T} {A_0} \qquad
A_- = \mat {a+2-2\cos\theta} {v} {v^T} {A_0},
$$ 
for $\rho=e^{i \theta}$ and some row vector $v$.
Using Lemma \ref{lem.tbordered} the result follows.
\end{proof}

\begin{lemma}
\lbl{lem.tbordered}
If $(A_+, A_-, A_0)$ is a $\rho$-bordered triple, and 
$\det(A_+) \det(A_-) \neq 0$, then
\begin{equation*}
\s(A_-)-\s(A_+)=
\begin{cases}
2 & \text{if} \quad \det(A_+) \det(A_-) < 0 \\
0 & \text{if} \quad \det(A_+) \det(A_-) > 0 .
\end{cases}
\end{equation*}
\end{lemma}

\begin{proof}
This is well-known for $\rho=1$;  \cite{C} and also \cite[Lemma 3.1]{Ga1}. 
We give a proof here for all $\rho$.

By a similarity transformation (that is a replacement of $B$ by
$P^\star B P$ where $P$ is an invertible matrix, and $P^\star$ is the 
conjugate transpose of $P$), we can assume that
$$
A_+= \mat {a} {v} {v^T} {0} \oplus D \qquad
A_-= \mat {a+2-2\cos\theta} {v} {v^T} {0} \oplus D, \qquad 
A_0= [0]^r \oplus D ,
$$ 
where $D$ is a nonsingular diagonal matrix, $[0]^r$ is the zero $r \times
r$ matrix,
$v$ is a $1\times r$ vector and $a$ a real number.

Since  the nullity (that is, the dimension of the kernel) and the
signature of the matrix  $\mat {b} {v} {v^T} {0}$ are given by: 
\begin{center}
\begin{tabular}{|l|c|c|c|}
\hline
        & $v=b=0$ & $v=0, b \neq 0$ & $ v \neq 0$ \\ \hline
nullity & $r+1$      &     $r$            &  $r-1$         \\ \hline
signature &  $0$     &   $\sgn(a)$        &   $0$           \\ \hline
\end{tabular}
\end{center}
\noindent
the result follows by a case-by-case argument. 
\end{proof}

The next theorem computes the jump function of a simple knot 
in terms of a relative sign of Alexander polynomials.
First, a preliminary definition. 

\begin{definition}
\lbl{def.sign}
If $f(x)$ is a real-valued analytic function of $x$ in a neighborhood of $a$,
we define {\em the sign of $f$ at $a$} $\sgn(f,a)$
to be the sign of the first nonvanishing Taylor series coefficient
(around $a$), if there
is such, and zero otherwise. In other words, if $f \neq 0$, we have:
$$
\sgn(f,a)=\sgn(f^{(n)}(a)) \, \in \{-1,1\},
$$
where $f^{(k)}(a)=0$ for $k < n$ and $f^{(n)}(a) \neq 0$.
\end{definition}

\begin{remark}
\lbl{rem.signf}
Notice that if $f(a) \neq 0$, then $\sgn(f,a)=\sgn(f(a))$, and that
if $a$ is a simple root, then $\sgn(f,a)=\sgn(f(a+\delta))=-\sgn(f(a-\delta))$
where $\delta$ is sufficiently small and positive.
\end{remark}

Fix a simple knot $K$ and a complex number $\rho=e^{i \theta} \in \DivD$.
Choose a planar projection of $K$ 
and a crossing (positive or negative). Then, $K=K^{\e}$, where
$\e \in \{+,-\}$ is the sign of the chosen crossing. 
Suppose that $\D(K^{-\e})(\rho) \neq 0$.
Such a projection and choice of crossing will be called $(\rho,K)$-{\em good}.

\begin{theorem}
\lbl{thm.jump}
Fix $(\rho,K)$ as above. For every $(\rho,K)$-good projection, we have
$$
\jump_{\rho}(K)=2 \e \, \sgn(\D(K^+), \theta) \sgn(\D(K^-),\theta) \in 
\{-2,2\}.
$$
\end{theorem}

\begin{proof}
Without loss of generality, let us assume $K=K^-$, that is $\e=-1$.
We will apply Lemma \ref{lem.skeins} twice to $\rho'=e^{i (\theta+\delta)}$
and $\rho''=e^{i (\theta-\delta)}$ for sufficiently small positive $\delta$.

Under these assumptions, we have that $\D(K^-)(\rho') \neq 0$ (since $\rho$ is
an isolated root of a polynomial) and $\D(K^+)(\rho') \neq 0$ 
(since $\D(K^+)(\rho)\neq 0$ by assumption), and similarly for $\rho''$.
Thus, the hypothesis of Lemma \ref{lem.skeins} are satisfied.
Applying Lemma \ref{lem.skeins} twice, we get
$$
\s_{\rho'}(K^-)-\s_{\rho'}(K^+)=
\begin{cases}
2 & \text{if} \quad \D(K^+)(\rho') \D(K^-)(\rho') < 0 \\
0 & \text{if} \quad \D(K^+)(\rho') \D(K^-)(\rho') > 0 
\end{cases}
$$
and 
$$
\s_{\rho''}(K^-)-\s_{\rho''}(K^+)=
\begin{cases}
2 & \text{if} \quad \D(K^+)(\rho'') \D(K^-)(\rho'') < 0 \\
0 & \text{if} \quad \D(K^+)(\rho'') \D(K^-)(\rho'') > 0 
\end{cases}
$$
Now, subtract and remember that $\s(K^+)$ is continuous at $\rho$ since
$\D(K^+)(\rho) \neq 0$. We get
$$
\jump_{\rho}(K^-)=
\begin{cases}
2 & \text{if} \quad \D(K^+)(\rho') \D(K^-)(\rho') < 0 \\
0 & \text{if} \quad \D(K^+)(\rho') \D(K^-)(\rho') > 0 
\end{cases}
-
\begin{cases}
2 & \text{if} \quad \D(K^+)(\rho'') \D(K^-)(\rho'') < 0 \\
0 & \text{if} \quad \D(K^+)(\rho'') \D(K^-)(\rho'') > 0 
\end{cases}
$$
Since $K$ is simple, it follows that $\D(K^-)(\rho') \D(K^-)(\rho'') <0$,
thus the cases $2-2$ or $0-0$ do not occur above.
Thus,
$$
\jump_{\rho}(K^-)=
\begin{cases}
2 & \text{if} \quad \D(K^+)(\rho') \D(K^-)(\rho') < 0 \\
-2 & \text{if} \quad \D(K^+)(\rho'') \D(K^-)(\rho'') < 0 
\end{cases}
$$
The result follows using Remark \ref{rem.signf}. Indeed,
$\sgn(\D(K^+),\rho)=\sgn(\D(K^+)(\rho))$ and
$\sgn(\D(K^-),\rho)=\sgn(\D(K^-)(\rho''))=-\sgn(\D(K^-)(\rho')$.
\end{proof}

\begin{theorem}
\lbl{thm.unique}
There is a unique invariant $j$ defined for a simple knot $K$ and $\rho \in
\DivD$ such that for every 
$(\rho,K)$-good projection we have:
$$
j_{\rho}(K)=2 \e \, \sgn(\D(K^+), \theta) \sgn(\D(K^-),\theta).
$$
\end{theorem}

\begin{proof}
In view of Theorem \ref{thm.jump}, we need to prove that there is
at most one such invariant.

Fix a simple knot $K$ and a complex number $\rho=e^{i \theta} \in \DivD$.
We need to prove that there exists a $(\rho,K)$-good projection.

Start with any planar projection of $K$ and a crossing. If it is not good,
apply Reidemaster moves II, which Frohman-Klassen call {\em threading}
and improve it to be good, using the proof of \cite[Theorem 6.2]{FK}.
\end{proof}

Thus, Conjecture \ref{conj.1} is equivalent to the following:

\begin{conjecture}
\lbl{conj.2}
$\mathrm{(a)}$ For every simple knot $K$, and every $\rho=e^{i \theta} 
\in \DivD$, we have $P(K)(\rho) \neq 0$.  \\
$\mathrm{(b)}$ Moreover, for every $(\rho,K)$-good projection we have:
$$
\sgn(P(K),\theta)=\e \, 
\sgn(\D(K^+), \theta)
\sgn(\D(K^-),\theta).
$$
\end{conjecture}

\section{Evidence}
\lbl{sec.evidence}

\subsection{Torus knots}
\lbl{sub.torusknots}

In this Section we will prove Conjecture \ref{conj.1} for torus knots.
Let $T_{a,b}$ denote the $(a,b)$ {\em torus knot}, where $a,b$ are coprime 
natural numbers. For example, $T(2,3)$ is the right-hand trefoil.

The Alexander polynomial of torus knots is given by:

\begin{eqnarray*}
\D(T_{a,b})(t) &=& \frac{(t^{ab/2}-t^{-ab/2})(t^{1/2}-t^{-1/2})}{
(t^{a/2}-t^{-a/2})(t^{b/2}-t^{-b/2})}.
\end{eqnarray*}

The roots of $\D(T_{a,b})$ on the unit circle are $ab$ complex roots of
unity which are not $a$ or $b$ order roots of unity. They are all simple.
Using a useful parametrization of them, following Kearton \cite[Sec.13]{K2}, 
we obtain that 
$$
\mathrm{Roots}_{\D(T_{a,b})}=\{t(m,n):=e^{2 \pi i (m/a+n/b)} \, 
| 0 < m < a, \,\,\, 0 < n < b \}.
$$
Since the jump function satisfies $\jump_{\rho}(K)=-\jump_{\bar\rho}(K)$,
we need only compute the jump at the points $t(m,n)$ where
$0 < m < a$, $0 < n < b$ and $m/a+n/b < 1$. In \cite[p.177]{K2} Kearton
computes the jump function of torus knots by
$$
\jump_{t_{m,n}}(T_{a,b})=
\begin{cases}
2 & \text{if} \quad m/a+n/b <\frac{1}{2} \\
-2 & \text{if} \quad \frac{1}{2} < m/a+n/b < 1.
\end{cases}
$$
In other words, we have:
$$
\jump_{\rho}(T_{a,b})=
\begin{cases}
-2 & \text{if} \quad \mathrm{Im}(\rho) > 0 \\
2 & \text{if} \quad \mathrm{Im}(\rho) < 0.
\end{cases}
$$

Now we discuss the $Q$ function of torus knots, which was originally
computed by Rozansky (see \cite[Eqn.(2.2)]{R1}), and most recently,
it has been recomputed by March\'e and Ohtsuki;
see \cite{Mr,Oh}. We understand that Bar-Natan has unpublished computations
of the Euler expansion of the Kontsevich integral of torus knots.

According to \cite[Eqn.(2.2)]{R1}, the $Q$ function of torus
knots is given by:

\begin{eqnarray*}
Q(T_{a,b})(t) &=& \frac{1}{4}\left( ab-\frac{a}{b}-\frac{b}{a} \right)
+ \frac{1}{ab} \frac{\D(T_{a,b})(t)}{ (t^{1/2}-t^{-1/2})}
\frac{\pt^2}{\pt x^2}\Big|_{x=0} \frac{t^{1/2}e^{x/2}-t^{-1/2}e^{-x/2}}{
\D(T_{a,b})(te^x)}
\end{eqnarray*}

Given an analytic function $f(t)$ let us define 
$$
g(t) =\frac{f(t)}{ (t^{1/2}-t^{-1/2})}
\frac{\pt^2}{\pt x^2}\Big|_{x=0} \frac{t^{1/2}e^{x/2}-t^{-1/2}e^{-x/2}}{
f(te^x)}
$$
We have that
$$
g(t) =
\frac{1}{8(t^{1/2}-t^{-1/2})} 
\frac{t f(t)^2 -f(t)^2 + 4 t^2 f(t) f''(t) -8 t^2 f(t) f'(t) 
-4t^3 f(t) f''(t) -8 t^2 (f'(t))^2 + 8 t^3 (f'(t))^2}{t^{1/2} f(t)^2}
$$
When we expand $g(e^{i \theta})$ around a root $\rho=e^{i \theta_0}$,
only the last two terms of the numerator 
contribute to the coefficient of $(\theta-\theta_0)^2$. That is,
$$
\coeff(g(e^{i \theta}), (\theta-\theta_0)^2)=
\frac{1}{8(t^{1/2}-t^{-1/2})} 
\frac{-8 t^2 (f'(t))^2 + 8 t^3 (f'(t))^2}{t^{1/2} f(t)^2}
\big|_{t=e^{i \theta_0}} =
 \frac{t^2 (f'(t))^2}{f(t)^2} \big|_{t=e^{i \theta_0}}.
$$
Now, suppose that $f(t)$ is a Laurrent polynomial with real coefficients
that satisfies $f(t)=f(t^{-1})$. Then, $f(t)=\sum_k a_k (t^k + t^{-k})$.
Thus, 
$$
t^2 (f'(t))^2=t^2 \left(\sum_k k a_k (t^{k-1}-t^{-k-1})\right)^2=
\left(\sum_k k a_k (t^k -t^{-k})\right)^2
$$
and if we substitute $t=e^{i \theta_0}$, we get
$$
t^2 (f'(t))^2 \big|_{t=e^{i \theta_0}}=-4\left(\sum_k k a_k \sin k 
\theta\right)^2
\leq 0.
$$
If $\theta_0$ is a simple root of $f(e^{i \theta})$ on the unit circle
(as is the case for the Alexander polynomial of torus knots), then the above
real number is negative.

This proves that
$$
\jjump_{\rho}(T_{a,b})=
\begin{cases}
-2 & \text{if} \quad \mathrm{Im}(\rho) > 0 \\
2 & \text{if} \quad \mathrm{Im}(\rho) < 0 
\end{cases}
$$
and confirms Conjecture \ref{conj.1} for torus knots.

\subsection{Operations on knots that preserve Conjecture \ref{conj.1}}
\lbl{sub.operations}

Let $f$ denote either the $Q$ function or the \sig\ function of a knot.
The following list describes some well-known properties of $f$.

\begin{itemize}
\item
If $-K$ denote the knot $K$ with opposite orientation, then $f(-K)=f(K)$.
\item
If $K^!$ denote the mirror image of $K$, then $f(K^!)=-f(K)$.
\item
If $\sharp$ denotes the connected sum of knots, then
$f(K_1 \sharp K_2) = f(K_1) + f(K_2)$.
\item
If $K^{(n)}$ denote the $(n,1)$ parallel of a knot $K$ with zero framing,
then $f(K^{(n)})(t)=f(K)(t^n)$. 
\end{itemize}

The stated behavior of the signature function under $(n,1)$ parallel
was proven by Kearton \cite{K1}, 
and for the $Q$ function was proven by Ohtsuki \cite[Prop. 3.1]{Oh}.

From this, it follows that if Conjecture \ref{conj.1} 
is true for a simple knot $K$,
then it is true for $-K$, $K^!$, $K^{(n)}$ (for all $n$). Furthermore, if
$K_1 \sharp K_2$ is simple, and Conjecture \ref{conj.1} is true for
$K_1$ and $K_2$, then it is also true for $K_1 \sharp K_2$.

\subsection{Knots with at most 8 crossings}
\lbl{sub.eight}

In this section we will verify Conjecture \ref{conj.1}
by computer calculations.

Rozansky has written a {\tt Maple} program that computes the $Q$ function
of a knot; see \cite{R2}.
We will use a minor modification {\tt Qfunction.mws} of Rozansky's 
program, adopted for our needs.

In {\tt Qfunction.mws}, the knot is described by a braid word.
For example, $[-1,3,3,3,2,1,1,-3,2]$ represents the braid 
$\s_1^{-1}\s_3^3 \s_2 \s_1^2 \s_3^{-1} \s_2$ whose closure is the $7_2$
knot in classical notation. The command $br1([-1,3,3,3,2,1,1,-3,2])$
gives a list whose first, second and third entries are the braid word, 
the polynomials $P(K)$ and $\D(K)$, where $z=t^{1/2}-t^{-1/2}$.
A sample output of the program is: 

{\small
\begin{verbatim}
> # the right trefoil 3_1
> br1([1,1,1]);
> 

                                     2     2    4
                    [[1, 1, 1], 1 + z , 2 z  + z ]

> # the 4_1 knot
> br1([1,-2,1,-2]);
> 

                                           2
                     [[1, -2, 1, -2], 1 - z , 0]
> # the 7_2 knot
> br1([-1,3,3,3,2,1,1,-3,2]);
> 

                                              2      2       4
       [[-1, 3, 3, 3, 2, 1, 1, -3, 2], 1 + 3 z , 12 z  + 14 z ]
> # 7_3
> br1([1,1,2,-1,2,2,2,2]);
> 

                                     2      4
  [[1, 1, 2, -1, 2, 2, 2, 2], 1 + 5 z  + 2 z ,

            2       4       6      8
        22 z  + 65 z  + 46 z  + 9 z ]
\end{verbatim}
}
For example, for the right hand trefoil, we have:

\begin{eqnarray*}
\D(K)&=& 1+z^2=t+\bt-1 \\
P(K) &=& 2z^2+z^4=t^2-2t+2-2\bt+\bt^2 \\
Q(K) &=& \frac{2z^2+z^4}{(1+z^2)^2}=\frac{t^2-2t+2-2\bt+\bt^2}{(t+\bt-1)^2}.
\end{eqnarray*}

The {\tt Mathematica program} {\tt JJump.m} computes the $\jjump$ function.
For example, we may launch the {\tt JJump.m} program from a Mathematica
session.

{\small
\begin{verbatim}
(math100)/home/stavros: math
Mathematica 5.0 for Sun Solaris (UltraSPARC)
Copyright 1988-2003 Wolfram Research, Inc.
 -- Motif graphics initialized --

In[1]:= << JJump.m

In[2]:= Poles[1+z^2,2z^2+z^4]

Solve::ifun: Inverse functions are being used by Solve, so some solutions may
     not be found; use Reduce for complete solution information.

Out[2]= {{0.16666666666666666667, -0.00844343197019481429}}
\end{verbatim}
}

We learn that the coefficient of $(\theta-\theta_0)^{-2}$ of 
$Q(3_1)(e^{2 \pi i \theta})$ (where $3_1$ is the right trefoil) 
around the root $\theta_0=0.1666666667$, 
is $-0.00844343197019481429$. This computes that $\jjump_{e^{2 \pi i\theta_0}}
(3_1)=-2$,
as needed.

Similarly, 

{\small
\begin{verbatim}
In[4]:= Poles[1+5z^2+2z^4,22z^2+65z^4+46z^6+9z^8]

Solve::ifun: Inverse functions are being used by Solve, so some solutions may
     not be found; use Reduce for complete solution information.

Out[4]= {{0.075216475230034463796, -0.00388836700144941422},

>    {0.27241752919082620707, -0.00542424178920663095}}
\end{verbatim}
}

We learn that the coefficient of $(\theta-\theta_i)^{-2}$ of  
$Q(7_3)(e^{2 \pi i \theta})$ 
around the roots $\theta_0=0.0752164$ and 
$\theta_1=0.27241752$ are $-0.003888367$ and $-0.0054242417$ respectively.
This computes the jump function  $\jjump_{e^{2 \pi i\theta_j}}(7_3)=-2$ for 
$j=0,1$.

Now, let us compute the jump function of a knot. In \cite{Or} Orevkov gives
a {\tt Mathematica} program {\tt sm.mat} which takes as input a braid
presentation of a knot, and gives as output a Seifert surface of a knot.
Launching the {\tt Jump.m} version of it in a Mathematica session produces

{\small
\begin{verbatim}
(math100)/home/stavros: math
Mathematica 5.0 for Sun Solaris (UltraSPARC)
Copyright 1988-2003 Wolfram Research, Inc.
 -- Motif graphics initialized --

In[1]:= << Jump.m

In[2]:= Jump[{1,1,1}]

InverseFunction::ifun:
   Inverse functions are being used. Values may be lost for multivalued
    inverses.

Solve::ifun: Inverse functions are being used by Solve, so some solutions may
     not be found; use Reduce for complete solution information.

Out[2]= {-2}
\end{verbatim}
}
which computes the jump function on the upper semicircle for the right trefoil
$3_1$.

{\small
\begin{verbatim}
In[3]:= Jump[{1,1,2,-1,2,2,2,2}]

InverseFunction::ifun:
   Inverse functions are being used. Values may be lost for multivalued
    inverses.

Solve::ifun: Inverse functions are being used by Solve, so some solutions may
     not be found; use Reduce for complete solution information.

Out[3]= {-2, -2}
\end{verbatim}
}
which computes the jump function on the upper semicircle for the $7_3$ knot.

This confirms the conjecture for the $3_1$ and $7_3$ knots.

In the appendix, We give the source code of two Mathematica programs, 
{\tt Jump.m} and {\tt JJump.m} which compute the $\jump$ and the $\jjump$ 
function of knots.

\appendix

\section{The {\tt JJump.m} program}
\lbl{sec.JJump}

{\small
\begin{verbatim}
     (* Poles[P,AP] computes the poles of the rational functions P/AP^2     *)
     (* at the roots of AP=0 on the unit circle. P,AP are polynomials in z  *)
     (* Poles2[P,AP] lists the coefficients of the Taylor expansion at      *)
     (* (t-a)^{-2}.                                                         *)
     (* Poles[P,AP] lists {a,coefficient of Taylor expansion at (t-a)^{-2}} *)
 

FF[x_]:=x[[2]];

Poles[AP_,P_]:=Module[
  {quotient,APt,roots,poles,k},
  quotient=Simplify[P/AP^2  /. (z->z^{1/2})  /. (z->2 Cos[2*Pi*t]-2 )];
  APt= Simplify[AP /. (z->z^{1/2})  /. (z->2 Cos[2*Pi*t]-2 )];
  roots=Select[Map[FF, Flatten[
      NSolve[APt == 0, t, 20]] ], 
	      1/2 > # > 0 &]; 
  poles={}; 
  Table[Flatten[{roots[[k]], Coefficient[Series[quotient,{t,roots[[k]],0}],
      t-roots[[k]],-2]}], {k,Length[roots]}]
]

     (* For the 3_1 knot:          Poles[1+z^2,2z^2+z^4]                    *)
     (* For the 4_1 knot:          Poles[1-z^2,0]                           *)
     (* For the 7_2 knot:          Poles[1+3z^2,12z^2+14z^4]                *)
     (* For the 7_3 knot:          Poles[1+5z^2+2z^4,22z^2+65z^4+46z^6+9z^8] *)
\end{verbatim}
}

\section{The {\tt Jump.m} program}
\lbl{sec.Jump}

{\small
\begin{verbatim}
    (* Computing the signature and jump function of knots presented as   *)
    (* closures of braids.                                               *)
    (* The signature of the right trefoil is SignatureBraid[{1,1,1}]=-2  *)
    (* SignatureM[A] of a matrix A is the signature of A+A^*             *)
    (* Jump[{1,1,1}] is the jumps of the signature of the right trefoil  *)

<< LinearAlgebra`MatrixManipulation`

<< sm.mat;

SignatureM[A_]:=Module[
  {eigen},
  eigen=Eigenvalues[N[A+ Transpose[Conjugate @ A],20]];
  Count[Sign @ eigen, 1]-Count[Sign @ eigen, -1]
]

SignatureBraid[brd_]:=Module[
  {m,V,eigen},
  m=Max[Abs @ brd]+1;
  V=N[SeifertMatrix[m,brd],20];
  SignatureM[V]
] 

FF[x_]:=x[[2]];

Jump[brd_]:=Module[
  {m,V,APs,hermitian,roots,k},
  m=Max[Abs @ brd]+1;
  V=N[SeifertMatrix[m,brd]];
  hermitian=(1-Exp[2*Pi*I*s])V+(1-Exp[-2*Pi*I*s]) Transpose[V];
  APs=N[Det[(Cos[2*Pi*s/2]+I Sin[2*Pi*s/2])V-(Cos[2*Pi*s/2]-I 
      Sin[2*Pi*s/2]) Transpose[V]],20];
  roots=Select[Map[FF, Flatten[
      NSolve[{APs == 0, Im[s]==0}, s, 15]] ], 1/2 > # > 0 &]; 
  If[Length[roots]==0, {}, Flatten[Table[SignatureM[hermitian /. 
      s->(roots[[k]]+1/1000) ] -SignatureM[hermitian /. 
s->(roots[[k]]-1/1000) ], {k,Length[roots]}]]]
]

     (* 7_3 knot     SignatureBraid[{1,1,2,-1,2,2,2,2}]    *)
     (* 7_5 knot     SignatureBraid[{1,1,1,1,2,-1,2,2}]    *)
     (* 8_2 knot     SignatureBraid[{-1,2,2,2,2,2,-1,2}]   *)
     (* 8_5 knot     SignatureBraid[{1,1,1,-2,1,1,1,-2}]   *)
     (* 8_15 knot    SignatureBraid[{1,1,-2,1,3,3,2,2,3}]  *)
     (* 7_3, 7_5, 8_2, 8_5, 8_15 have signature       -4   *)

\end{verbatim}
}

\ifx\undefined\bysame
	\newcommand{\bysame}{\leavevmode\hbox
to3em{\hrulefill}\,}
\fi


\begin{thebibliography}{[EMSS]}

\bibitem[AM]{AM} S. Akbulut, J. C.  McCarthy, 
        {\em Casson's invariant for oriented homology 3-spheres: 
        an exposition},  
        Princeton Math Notes, Princeton, 1990.

\bibitem[BG]{BG} D. Bar-Natan, S. Garoufalidis,
        {\em On the Melvin-Morton-Rozansky conjecture}, 
        Inventiones, {\bf 125} (1996) 103--133.

\bibitem[BLT]{BLT} \bysame, T.T.Q. Le and D. Thurston,
        {\em Two applications of elementary knot theory to Lie algebras
        and Vassiliev invariants}, 
         Geom. Topol.  {\bf 7}  (2003) 1--31.

\bibitem[C]{C} J. Conway,
        {\em An enumeration of knots and links and some of their algebraic
        properties},
        Computational problems in abstract algebra, Pergamon Press,
        New-York 1970, 329--358.

\bibitem[Ga1]{Ga1} S. Garoufalidis,
        {\em Signatures of links and finite type invariants of cyclic
        branched covers}, 
        Contemporary Math. {\bf 231} (1999) 87--97.

\bibitem[GR]{GR} \bysame and L. Rozansky, 
        {\em The loop expansion of the Kontsevich integral, 
        abelian invariants of knots and $S$-equivalence}, 
        preprint 2000, {\tt math.GT/0003187}, to appear in Topology.

\bibitem[GK1]{GK1} \bysame and A. Kricker,
        {\em A rational noncommutative invariant of boundary links}, 
        preprint 2001, {\tt math.GT/0105028}.

\bibitem[GK2]{GK2} \bysame and \bysame,
        {\em Finite type invariants of cyclic branched covers}, 
        preprint 2001, {\tt math.GT/0107220}.

\bibitem[GL]{GL} \bysame and TTQ. Le,
        {\em The colored Jones function is $q$-holonomic}
        preprint 2003, {\tt math.GT/0309214}.

\bibitem[Ga2]{Ga2} \bysame,
        Programs {\tt Qfunction.mws}, {\tt Jump.m} and {\tt JJump.m},
        available upon request.

\bibitem[Ga3]{Ga3} \bysame,
        {\em On the characteristic and deformation varieties of a knot},
        preprint 2003 {\tt math.GT/0306230}.

\bibitem[FK]{FK} C. Frohman and E. Klassen,
        {\em  Deforming representations of knot groups in ${\rm SU}(2)$},
        Comment. Math. Helv. {\bf 66} (1991) 340--361.

\bibitem[H1]{H1} C. Herald,
        {\em Existence of irreducible representations of knot complements
        with nonconstant equivariant signature},
        Math. Annalen {\bf 309} (1997) 21--35.

\bibitem[H2]{H2} \bysame,
        {\em Flat connections, the Alexander invariant and Casson's invariant},
        Comm. Anal. Geom. {\bf 5} (1997) 93--120.

\bibitem[K1]{K1} C. Kearton,
        {\em The Milnor signatures of compound knots},
        Proc. Amer. Math. Soc. {\bf 76} (1979) 157--160.

\bibitem[K2]{K2} \bysame,
        {\em Signatures of knots and the free differential calculus},
        Quart. J. Math. Oxford Ser. {\bf 30}  (1979) 157--182.

\bibitem[L]{L} J. Levine,
        {\em Invariants of knot cobordism},
        Inventiones Math. {\bf 8} (1969) 98--110.

\bibitem[Li]{Li} X-S. Lin,
        {\em A knot invariant via representation spaces},
        J. Differential Geom.  {\bf 35}  (1992)  337--357.

\bibitem[Mr]{Mr} J. March\'e,
        {\em On Kontsevich integral of torus knots},
        preprint 2003 {tt math.GT/0310111}.

\bibitem[Ma]{Ma} B. Mazur, 
        {\em The theme of $p$-adic variation},   
        Mathematics: frontiers and perspectives,  433--459, 
        Amer. Math. Soc., Providence, RI, 2000.

\bibitem[M1]{M1} J. Milnor, 
        {\em Infinite cyclic coverings}, 
        Conference on the Topology of Manifolds, 
        Michigan State University (1967)  115--133.

\bibitem[M2]{M2} \bysame,
        {\em On isometries of inner product spaces},
        Invent. Math. {\bf 8} (1969) 83--97.

\bibitem[Oh]{Oh} T. Ohtsuki,
        {\em A cabling formula for the 2-loop polynomial of knots},
        preprint 2003.

\bibitem[Or]{Or} S. Orevkov,
        {\em Classification of flexible $M$-curves of degree 8 up to isotopy},
        Geom. Funct. Anal. {\bf 12}  (2002) 723--755. 

\bibitem[R1]{R1} L. Rozansky,
        {\em  Higher order terms in the Melvin-Morton expansion of the 
        colored Jones polynomial},  
        Comm. Math. Phys. {\bf 183} (1997) 291--306.

\bibitem[R2]{R2} L. Rozansky,
        {\em The universal $R$-matrix, Burau Representation and the 
        Melvin-Morton expansion of the colored Jones polynomial},
        Adv. Math. {\bf 134} (1998) 1--31.

\bibitem[R3]{R3} \bysame,
        Computer programs {\tt pol1.mws, pol2.mws} in Maple code, available
        at \newline {\tt http://www.math.yale.edu/\~{}rozansky} 

\end{thebibliography}
\end{document}